\documentclass[psamsfonts]{amsart}

\usepackage{amssymb,amsfonts}
\usepackage[all,arc]{xy}
\usepackage{enumerate}
\usepackage{mathrsfs}
\usepackage{hyperref}
\usepackage{mathtools}
\usepackage{comment}
\usepackage{esint}
\hypersetup{
	pdftitle={The Maximum Principle and Optimal Gradient Estimates of Solutions to the Insulated Conductivity Problem},
	pdfsubject={Mathematics, Differential Calculus, Derivatives, Differential equations},
	pdfauthor={Saibo Chen},
	pdfkeywords={Maximum Principle, Insulated Conductivity Problem, Bispherical Coordinates, Gradient Estimates}
}
\hypersetup{pdfborder={0 0 0}}
\usepackage{bm}
\usepackage{pgfplots}

\newtheorem{thm}{Theorem}[section]

\newtheorem{prop}{Proposition}[section]
\newtheorem{lem}{Lemma}[section]
\newtheorem{rk}{Remark}[section]

\theoremstyle{definition}

\theoremstyle{remark}

\makeatletter
\let\c@equation\c@thm
\makeatother
\numberwithin{equation}{section}

\title[The Maximum Principle and the Insulated Conductivity Problem]{The Maximum Principle and Optimal Gradient Estimates of Solutions to the Insulated Conductivity Problem}

\author{Saibo Chen}
\address{School of Mathematical Sciences, Fudan University, Shanghai 200433, China}
\email{sbchen25@m.fudan.edu.cn}

\begin{document}
	
	\begin{abstract}
		We study the insulated conductivity problem with two unit spherical inclusions embedded in a bounded domain in $\mathbb{R}^n$ with $n\ge3$. The gradient of solutions may blow up as the distance between the two inclusions goes to zero. We introduce bispherical coordinates and expand the solution in spherical harmonics. Then we use the maximum principle to give a new proof of the optimal gradient estimate obtained by Dong-Li-Yang.
	\end{abstract}

	\maketitle
	\section{Introduction}
		Let $B_1$ and $B_2$ be two balls of radius $1$, compactly contained in a bounded open set $\Omega\subset\mathbb{R}^n$ with smooth boundary. Write $x=(x',x_n)\in\mathbb{R}^n$, where $x'=(x_1,...,x_{n-1})$. Assume that $B_i$ is centered at $(0',(-1)^{i+1}(\varepsilon+1))$, where $\varepsilon$ is a small parameter. We further assume that the distance from $B_1\cup B_2$ to $\partial\Omega$ is bounded below by a positive constant independent of $\varepsilon$.
		
		For a fixed constant $\kappa\in(0,\infty)$, let
		\begin{equation}\label{conductivity}
			\tilde{a}(x)\coloneqq
			\begin{cases}
				\kappa\quad{\rm in}\ B_1\cup B_2,\\
				1\quad{\rm in}\ \Omega\setminus\overline{B_1\cup B_2}
			\end{cases}
		\end{equation}
		be the conductivity, and let $u(x)$ be the voltage potential. Then the conductivity problem is given by
		\begin{equation}\label{general conductivity}
			\begin{cases}
				\nabla\cdot(\tilde{a}\nabla u)=0\quad{\rm in}\ \Omega,\\
				u=\varphi\quad{\rm on}\ \partial\Omega,
			\end{cases}
		\end{equation}
		where $\varphi$ is a smooth function on $\partial\Omega$. We want to characterize the asymptotic singular behavior of $|\nabla u|$ as $\varepsilon\to0$.
		
		Li and Vogelius \cite{generalk} proved that $|\nabla u|$ remains bounded independently of $\varepsilon$ for a fixed $\kappa$, and this was extended to elliptic systems by Li and Nirenberg \cite{generalsys}. If we consider the limiting cases $\kappa\to\infty$ (perfect conductors) or $\kappa\to0$ (insulators), the results become more complicated. We refer the reader to the references \cite{Ammariqqb,Ammari2,Ammari1,perfect-alldim,insulated-dim2,green function,insulated-alldim general,insulated-alldim,KLYP,insulated-Harnack,perfect-dim2-general,Yun another,Yun Ball}.
		
		We point out that Weinkove \cite{maxprin} developed a maximum principle method to get effective (not optimal) estimates for the insulated conductivity problem when $n\ge4$. Li and Weinkove \cite{maximum principle with general k} used a similar method to reprove the result in \cite{generalk} in two dimensions with circular inclusions. Moreover, they speculated that the maximum principle alone can be used to reprove all the existing results we listed above. Motivated by this idea, we give a new proof of the optimal gradient estimate for the insulated conductivity problem with two spherical inclusions of equal radius in dimensions $n\ge3$. We first use bispherical coordinates, where the boundary condition becomes simpler. Then we take a decomposition of the solution in spherical harmonics. Finally, the maximum principle method is used to estimate each term.
		
		We now introduce our notation and state the main theorem. The insulated conductivity problem is
		\begin{equation}\label{insulated conductivity}
			\begin{cases}
				\Delta u=0\quad{\rm in}\ \Omega\setminus\overline{B_1\cup B_2},\\
				\partial_{\nu}u=0\quad{\rm on}\ \partial B_1\cup\partial B_2,\\
				u=\varphi\quad{\rm on}\ \partial\Omega,
			\end{cases}
		\end{equation}
		where $\nu$ is the unit normal pointing into $B_1$ or $B_2$. Choose a constant $r_0>0$ which is sufficiently small and independent of $\varepsilon$. Then choose a small $\varepsilon_0>0$ which may depend on $r_0$. Throughout this paper, we assume $\varepsilon\in(0,\varepsilon_0)$. We may decrease $r_0$ and $\varepsilon_0$ whenever necessary. A positive constant is called universal if it is independent of $\varepsilon,k,i$, and $\delta$ (defined below) but may depend on $\Omega,\varphi,n,r_0$ and $\varepsilon_0$. The expression $X\lesssim Y$ implies that there exists a universal $C$ such that $X\le CY$, and $X\sim Y$ implies that both $X\lesssim Y$ and $Y\lesssim X$ hold. The notation $\nabla v$ refers to the Euclidean gradient of $v$.
		
		Let
		$$h(\rho)\coloneqq\varepsilon+1-\sqrt{1-\rho^2}.$$
		Define
		$$\Omega_0\coloneqq\{x\in\Omega\setminus\overline{B_1\cup B_2}:|x'|<r_0,\ |x_n|<h(|x'|)\}.$$
		By the standard elliptic estimates, we have
		\begin{lem}
			Let $u$ be the solution of (\ref{insulated conductivity}). For every positive integer $t$, there exists a constant $C_t$ (independent of $\varepsilon$) such that
			$$|\nabla^tu|\le C_t\quad{\rm on}\ \overline{\Omega\setminus(\Omega_0\cup B_1\cup B_2)}.$$
		\end{lem}
		So it suffices to focus on the narrow region near the origin. Using the maximum principle method, we have
		\setcounter{thm}{0}
		\begin{thm}
			Assume $n\ge3$, and let $u$ be the solution of (\ref{insulated conductivity}). Then
			$$|\nabla u(x)|\lesssim(\varepsilon+|x'|^2)^{\frac{\alpha-1}{2}}\quad{\rm in}\ \Omega_0,$$
			where
			$$\alpha=\frac{-(n-1)+\sqrt{(n-1)^2+4(n-2)}}{2}\in(0,1).$$
		\end{thm}
		\begin{rk}
			For $n\ge2$, it can be proved by the maximum principle that
			$$|\nabla u|\lesssim(\varepsilon+|x'|^2)^{-1/2}\quad{\rm in}\ \Omega_0,$$
			see \cite{maxprin}. Moreover, this rate is optimal when $n=2$. So we only consider $n\ge3$ in this paper.
		\end{rk}
		\begin{rk}
			The blow-up rate in Theorem 1.1 is the same as in \cite{insulated-alldim}, which is proved to be optimal. We only consider the case for spherical inclusions of equal radius. In more general cases, our method may not work.
		\end{rk}
		The rest of this paper is organized as follows. In Section 2, we introduce the bispherical coordinates and decompose $u$. In Section 3, we construct some auxiliary functions. In Sections 4 and 5, we use these functions to give estimates for each term in the decomposition.

	\section{Bispherical coordinates and the decomposition of $u$}
		We first recall the cylindrical coordinates $(\rho,\omega,x_n)$, where
		$$\rho=|x'|,\quad\omega\in\mathbb{S}^{n-2}.$$
		In the $\rho O x_n$ half-plane, let $P_{\pm}=(0,\pm a)$, where
		\begin{equation}\label{a}
			a=\sqrt{\varepsilon(\varepsilon+2)}.
		\end{equation}
		For a point $P=(\rho,x_n)$, let
		$$\xi\coloneqq\log\frac{|PP_-|}{|PP_+|}\in(-\infty,+\infty),\quad\eta\coloneqq\angle P_-PP_+\in(0,\pi].$$
		This defines the bispherical coordinates $(\xi,\eta,\omega)$. One can check that
		\begin{equation}\label{xi eta}
			\rho=\frac{a\sin\eta}{\cosh\xi-\cos\eta},\quad x_n=\frac{a\sinh\xi}{\cosh\xi-\cos\eta}.
		\end{equation}
	
		Now we consider the shape of $\Omega_0$ in this coordinate system. A direct calculation shows
		\begin{equation}\label{xi=xi0 is B}
			\rho^2+(x_n-a\coth\xi)^2=a^2\sinh^{-2}\xi.
		\end{equation}
		Choose $\xi_0$ such that $\sinh\xi_0=a$. Then
		$$a\coth\xi_0=\varepsilon+1,\quad a^2\sinh^{-2}\xi_0=1.$$
		Hence, the hypersurfaces $\{\xi=\pm\xi_0\}$ represent $\partial B_1$ and $\partial B_2$.
		
		For $\eta=\pi$, we have
		$$\rho=0,\quad x_n=a\tanh\frac{\xi}{2}\in(-a,a).$$
		Hence, $\{\eta=\pi\}$ represents the open line segment $P_-P_+$.
		
		By (\ref{xi eta}),
		$$\rho=a\cot\frac{\eta}{2}\cdot\frac{1-\cos\eta}{\cosh\xi-\cos\eta}.$$
		Choose $\eta_0$ such that
		$$a\cot\frac{\eta_0}{2}=4r_0,$$
		where $r_0$ is the same as in the definition of $\Omega_0$. Then we can check
		\begin{equation}\label{rho eta}
			\frac{1}{2}\le\frac{1-\cos\eta}{\cosh\xi-\cos\eta}\le1\quad{\rm and}\quad\rho\sim a\cot\frac{\eta}{2}\quad{\rm for}\ |\xi|\le\xi_0,\ \eta\in[\eta_0,\pi].
		\end{equation}
		By the above argument, the domain corresponding to $\Omega_0$ in bispherical coordinates is contained in $R\times\mathbb{S}^{n-2}$, where
		$$R=\{(\xi,\eta):|\xi|<\xi_0,\eta_0<\eta\le\pi\}.$$
		Hence, we only need to study the problem in $R\times\mathbb{S}^{n-2}$. 
		
		Now we consider the gradient and Laplacian in the bispherical coordinates. We have
		\begin{equation}\label{bigradient}
			|\nabla u|^2=\frac{1}{\rho^2}\left[\sin^2\eta(u_{\xi}^2+u_{\eta}^2)+|\nabla_{\mathbb{S}^{n-2}}u|^2\right]
		\end{equation}
		and
		\begin{equation}\label{bilaplace}
			\begin{aligned}
				\Delta u=\frac{(\cosh\xi-\cos\eta)^2}{a^2}\biggl[&u_{\xi\xi}+u_{\eta\eta}-(n-2)\frac{\sinh\xi}{\cosh\xi-\cos\eta}u_{\xi}\\
				&+(n-2)\left(\cot\eta-\frac{\sin\eta}{\cosh\xi-\cos\eta}\right)u_{\eta}+\frac{\Delta_{\mathbb{S}^{n-2}}u}{\sin^2\eta}\biggr].
			\end{aligned}
		\end{equation}
		Let $\{Y_{k,i}\}_{k\ge0,\ 1\le i\le N(k)}$ be a real orthonormal basis of $L^2(\mathbb{S}^{n-2})$, where $Y_{k,i}$ is a $k$-th degree spherical harmonic, that is,
		\begin{equation}\label{Yki}
			-\Delta_{\mathbb{S}^{n-2}}Y_{k,i}=k(k+n-3)Y_{k,i}
		\end{equation}
		and $N(k)$ is the dimension of the corresponding space. We point out that when $k=0$, we have $N(k)=1$ and $Y_{0,1}$ is a constant. Therefore, we have a decomposition
		\begin{equation}\label{fourier series}
			u(\xi,\eta,\mathbf{\omega})=u_0(\xi,\eta)+\sum_{k=1}^{\infty}\sum_{i=1}^{N(k)}F_{k,i}(\xi,\eta)Y_{k,i}(\omega),
		\end{equation}
		where
		$$u_0=\fint_{\mathbb{S}^{n-2}}u\ {\rm d}\omega,\quad F_{k,i}=\int_{\mathbb{S}^{n-2}}uY_{k,i}\ {\rm d}\omega.$$
		For convenience, let
		$$u_{k,i}(\xi,\eta,\omega)\coloneqq F_{k,i}(\xi,\eta)Y_{k,i}(\omega).$$
		By (\ref{bigradient}) we have
		\begin{equation}\label{nabla k,i in bispherical coordinates}
			\begin{aligned}
				|\nabla u_{k,i}|\lesssim\frac{1}{\rho}\bigg[&\sin\eta\Vert Y_{k,i}\Vert_{L^{\infty}(\mathbb{S}^{n-2})}(|F_{k,i,\xi}|+|F_{k,i,\eta}|)\\
				&+\Vert\nabla_{\mathbb{S}^{n-2}}Y_{k,i}\Vert_{L^{\infty}(\mathbb{S}^{n-2})}|F_{k,i}|\bigg].
			\end{aligned}
		\end{equation}
		Hence, it suffices to estimate $|\nabla u_0|$, $|F_{k,i}|$, $|F_{k,i,\xi}|$ and $|F_{k,i,\eta}|$ ($k\ge1$). Sections 3 and 4 are devoted to estimating the quantities related to $F_{k,i}$, while Section 5 gives an estimate for $|\nabla u_0|$.

	\section{Construction and estimates of auxiliary functions}
		We first consider the estimates for $F_{k,i}$, $k\ge1$. For convenience, we denote
		$$\lambda_k\coloneqq k(k+n-3).$$
		Let
		\begin{equation}\label{Lk}
			\begin{aligned}
				L_{k}&\coloneqq\sin^2\eta(\partial^2_{\xi\xi}+\partial^2_{\eta\eta})-(n-2)\frac{\sin^2\eta\sinh\xi}{\cosh\xi-\cos\eta}\partial_{\xi}\\
				&+(n-2)\sin\eta\left(\cos\eta-\frac{\sin^2\eta}{\cosh\xi-\cos\eta}\right)\partial_{\eta}-\lambda_k.
			\end{aligned}
		\end{equation}
		By (\ref{bilaplace}), (\ref{fourier series}) and (\ref{Yki}) we have
		\begin{equation}\label{Fki}
			\begin{cases}
				L_{k}F_{k,i}=0\quad {\rm in}\ R,\\
				F_{k,i,\xi}=0,\quad |\xi|=\xi_0.
			\end{cases}
		\end{equation}
		Our strategy is to find a simple function $\psi_k$ so that $L_k\psi_k\le0$. 
		Then the maximum principle makes it possible to estimate $F_{k,i}$ and its derivatives.
		
		The operator $L_k$ can be rewritten as
		\begin{equation}
			\begin{aligned}
				L_{k}=&\sin^2\eta\partial^2_{\xi\xi}-(n-2)\frac{\sin^2\eta\sinh\xi}{\cosh\xi-\cos\eta}\partial_{\xi}+(n-2)\sin\eta\frac{(\cosh\xi-1)(\cos\eta+1)}{\cosh\xi-\cos\eta}\partial_{\eta}\\
				&+\left[\sin^2\eta\partial_{\eta\eta}-(n-2)\sin\eta\partial_{\eta}-\lambda_k\right].
			\end{aligned}\nonumber
		\end{equation}
		We rearrange the derivative terms with respect to $\eta$. The last term in the first line has a nonnegative coefficient, while the second line depends only on $\eta$. Before we give a construction of $\psi_k$, we introduce a lemma \cite[Theorem 4.5]{ode lemma}.
		\begin{lem}
			Consider the equation
			\begin{equation}\label{ode lemma}
				u''(x)+p(x)u'+q(x)u=0.
			\end{equation}
			Suppose the coefficients $p(x)$ and $q(x)$ have poles at $x=0$ of order at most one and two, respectively. Define
			$$p_0\coloneqq\lim_{x\to0}xp(x),\quad q_0=\lim_{x\to0}x^2q(x),$$
			and
			$$\beta_{1,2}=\frac{1-p_0\pm\sqrt{(p_0-1)^2-4q_0}}{2}.$$
			Assume ${\rm Re}\beta_1\ge{\rm Re}\beta_2$. Then we have
			\begin{enumerate}
				\item If $\beta_1-\beta_2\notin\mathbb{N}_0$, a fundamental system of solutions is given by
				$$u_j(x)=x^{\beta_j}h_j(x),\quad j=1,2,$$
				where the functions $h_j(x)$ are analytic near $x=0$ and satisfy $h_j(0)=1$.
				\item If $\beta_1-\beta_2=m\in\mathbb{N}_0$, a fundamental system of solutions is given by
				$$u_1(x)=x^{\beta_1}h_1(x),$$
				$$u_2(x)=x^{\beta_2}h_2(x)+c\log(x)u_1(x),$$
				where the functions $h_j(x)$ are analytic near $x=0$ and satisfy $h_j(0)=1$. The constant $c$ might be zero unless $m=0$.
			\end{enumerate}
		\end{lem}
		Using this lemma, we have
		\begin{lem}
			The equation
			\begin{equation}\label{equation of psi}
				\sin^2\eta\psi_k''-(n-2)\sin\eta\psi_k'-\lambda_k\psi_k=0,\quad\eta\in(0,\pi)
			\end{equation}
			has a solution satisfying $\psi_k>0$,\ $\psi_k(\eta_0)=1$, and $\psi_k'<0$.
		\end{lem}
		\begin{proof}
			Set
			$$\tilde{\eta}\coloneqq\pi-\eta,\quad\tilde{\psi}_k(\tilde{\eta})\coloneqq\psi_k(\pi-\tilde{\eta}).$$
			Then (\ref{equation of psi}) can be written as
			\begin{equation}\label{equation of tilde psi}
				\tilde{\psi}_k''+\frac{n-2}{\sin\tilde{\eta}}\tilde{\psi}_k'-\frac{\lambda_k}{\sin^2\tilde{\eta}}\tilde{\psi}_k=0,\quad\tilde{\eta}\in(0,\pi).
			\end{equation}
			We can check that the hypotheses of Lemma 3.1 (2) hold, so in a neighborhood of $0$, one solution of (\ref{equation of tilde psi}) can be given by
			$$\tilde{\psi}_k(\tilde{\eta})=\tilde{\eta}^k\tilde{h}_{k,0}(\tilde{\eta}).$$
			By the extension theorem for ODEs, we can extend this solution to $(0,\pi)$. Then we get a solution $\psi_k$ for (\ref{equation of psi}), and
			\begin{equation}\label{psik near pi}
				\psi_k(\eta)=(\pi-\eta)^kh_{k,0}(\eta)
			\end{equation}
			in a neighborhood of $\pi$ with
			$$h_{k,0}(\eta)=\tilde{h}_{k,0}(\pi-\eta).$$
			By the properties of $h_{k,0}$ in Lemma 3.1, we have $\psi_k>0$ and $\psi_k'<0$ in a neighborhood of $\pi$. Now we claim these hold in the whole interval $(0,\pi)$.
			
			In fact, if $\psi_k>0$ does not always hold, suppose the largest zero is $\eta_*\in(0,\pi)$. Then $\psi_k$ is positive in $(\eta_*,\pi)$ and $\psi_k'(\eta_*)\ge0$. However, by (\ref{equation of psi}) we have
			\begin{equation}\label{another equation of psi}
				\left(\cot^{n-2}\frac{\eta}{2}\psi_k'\right)'=\lambda_k\frac{\cot^{n-2}\frac{\eta}{2}}{\sin^2\eta}\psi_k.
			\end{equation}
			Integrating from $\eta_*$ to $\pi$ shows
			$$-\cot^{n-2}\frac{\eta_*}{2}\psi_k'(\eta_*)=\int_{\eta_*}^{\pi}\lambda_k\frac{\cot^{n-2}\frac{\eta}{2}}{\sin^2\eta}\psi_k\ {\rm d}\eta>0,$$
			which is a contradiction. Hence,
			$$\psi_k(\eta)>0,\quad\eta\in(0,\pi).$$
			Then integrating (\ref{another equation of psi}) from $\eta$ to $\pi$ implies
			$$\psi_k'(\eta)<0,\quad\eta\in(0,\pi).$$
			Finally, we can multiply $\psi_k$ by a positive constant so that $\psi_k(\eta_0)=1$. 
		\end{proof}
		Define
		$$\phi_k(\eta)\coloneqq-\sin\eta\frac{\psi_k'}{\psi_k}.$$
		In Section 4, we will use some estimates for this function.
		\begin{lem}
			We have
			\begin{equation}\label{estimate of phik}
				\alpha_k\le\phi_k\le k,\quad\eta\in[0,\pi],
			\end{equation}
			where
			\begin{equation}\label{alphak}
				\alpha_k=\frac{-(n-1)+\sqrt{(n-1)^2+4\lambda_k}}{2}.
			\end{equation}
			We point out that $\alpha_1=\alpha$, where $\alpha$ is the constant in Theorem 1.1.
		\end{lem}
		\begin{proof}
			We see from (\ref{psik near pi}) that
			$$\lim_{\eta\to\pi}\phi_k(\eta)=k.$$
			As for the case $\eta\to0$, we can check that the hypotheses of Lemma 3.1 (1) hold. So there exist constants $C_1,\ C_2$ such that $\psi_k$ can be written as
			\begin{equation}\label{psik near 0}
				\psi_k(\eta)=C_1\eta^{-\alpha_k}h_{k,1}(\eta)+C_2\eta^{\alpha_k+n-1}h_{k,2}(\eta)
			\end{equation}
			in a neighborhood of $0$. Since $\psi_k>0$ and $\psi_k'<0$, we have $C_1>0$. Then a direct calculation shows
			$$\lim_{\eta\to0}\phi_k(\eta)=\alpha_k.$$
			So $\phi_k$ extends continuously to $[0,\pi]$. By (\ref{equation of psi}), we have
			\begin{equation}\label{equation of phi}
				\sin\eta\phi_k'=\phi_k^2+(n-2+\cos\eta)\phi_k-\lambda_k.
			\end{equation}
			We first show $\phi_k\le k$. If there exists $\eta^*\in(0,\pi)$ such that
			$$\max_{\eta\in[0,\pi]}\phi_k=\phi_k(\eta^*)>k,$$
			then (\ref{equation of phi}) implies that
			$$0=\phi_k^2(\eta^*)+(n-2+\cos\eta^*)\phi_k(\eta^*)-\lambda_k>0,$$
			which is a contradiction. Hence, $\phi_k\le k$. We can prove $\phi_k\ge\alpha_k$ in the same way.
		\end{proof}
		The next lemma gives an estimate for $\psi_k$, which plays an important role.
		\begin{lem}
			Define
			\begin{equation}\label{Psik}
				\Psi_k(\eta)\coloneqq\sin^{\alpha_k}\frac{\eta_0}{2}\cos^{-k}\frac{\eta_0}{2}\sin^{-\alpha_k}\frac{\eta}{2}\cos^k\frac{\eta}{2}.
			\end{equation}
			Then
			\begin{equation}\label{psik and Psik}
				\psi_k\le\Psi_k,\quad\eta\in[\eta_0,\pi].
			\end{equation}
		\end{lem}
		\begin{proof}
			We have
			$$\Psi_k'=\frac{\Psi_k}{2}\left(-k\cos^{-1}\frac{\eta}{2}\sin\frac{\eta}{2}-\alpha_k\cos\frac{\eta}{2}\sin^{-1}\frac{\eta}{2}\right),$$
			$$\Psi_k''=\frac{\Psi_k}{4}\left[(k^2-k)\cos^{-2}\frac{\eta}{2}\sin^2\frac{\eta}{2}+(\alpha_k^2+\alpha_k)\cos^2\frac{\eta}{2}\sin^{-2}\frac{\eta}{2}+(2k\alpha_k+\alpha_k-k)\right].$$
			Hence,
			$$\sin^2\eta\Psi_k''-(n-2)\sin\eta\Psi_k'-k(k+n-3)\Psi_k=-(k-\alpha_k)^2\sin^2\frac{\eta}{2}\cos^2\frac{\eta}{2}\Psi_k.$$
			Therefore,
			\begin{equation}\label{psi-Psi}
				\sin^2\eta(\psi_k-\Psi_k)''-(n-2)\sin\eta(\psi_k-\Psi_k)'-\lambda_k(\psi_k-\Psi_k)\ge0.
			\end{equation}
			Suppose
			$$\max_{\eta\in[\eta_0,\pi]}(\psi_k-\Psi_k)=\psi_k(\eta_{**})-\Psi_k(\eta_{**})>0.$$
			Since
			$$(\psi_k-\Psi_k)(\eta_0)=(\psi_k-\Psi_k)(\pi)=0,$$
			we have
			$$\eta_{**}\in(\eta_0,\pi),\quad(\psi_k-\Psi_k)''(\eta_{**})\le0,\quad(\psi_k-\Psi_k)'(\eta_{**})=0.$$
			Substituting these into (\ref{psi-Psi}) leads to a contradiction. Hence, (\ref{psik and Psik}) holds.
		\end{proof}

	\section{Estimates of $|\nabla u_{k,i}|$}
		Define
		\begin{equation}
			\begin{aligned}
				a_{k,i}\coloneqq\sup_{|\xi|\le\xi_0,\eta=\eta_0}\bigg[\sin\eta_0\Vert Y_{k,i}\Vert_{L^{\infty}}(|F_{k,i,\xi}|+|F_{k,i,\eta}|)+\Vert \nabla_{\mathbb{S}^{n-2}}Y_{k,i}\Vert_{L^{\infty}}|F_{k,i}|\bigg]
			\end{aligned}\nonumber
		\end{equation}
		and
		$$E_{k,i}\coloneqq\sin^2\eta(F_{k,i,\xi}^2+F_{k,i,\eta}^2)+A_{k,i}F_{k,i}^2,$$
		where $A_{k,i}>0$ is a constant to be determined. In this section, we first use $a_{k,i}$ and $\psi_k$ to control $|F_{k,i}|$ and then estimate $E_{k,i}$.
		\subsection{Estimates of $|F_{k,i}|$}
			\begin{lem}
				$$|F_{k,i}|\le\frac{a_{k,i}}{\Vert\nabla_{\mathbb{S}^{n-2}}Y_{k,i}\Vert_{L^{\infty}}}\psi_k\quad{\rm in}\ R.$$
			\end{lem}
			\begin{proof}
				By the definition of $L_k$ and $\psi_k$, we have $L_k\psi_k\le0$. Hence, we see from (\ref{Fki}) that
				\begin{equation}\label{Cpsi-F}
					L_k\left(C_{k,i}\psi_k\pm F_{k,i}\right)\le0\quad{\rm in}\ R_{\delta},
				\end{equation}
				where
				$$C_{k,i}=\frac{a_{k,i}}{\Vert\nabla_{\mathbb{S}^{n-2}}Y_{k,i}\Vert_{L^{\infty}(\mathbb{S}^{n-2})}}$$
				and
				$$R_{\delta}=\{(\xi,\eta):|\xi|<\xi_0,\ \eta_0<\eta<\pi-\delta\}.$$
				Then $L_k$ is uniformly elliptic in $R_{\delta}$. By (\ref{Fki}), we have
				$$F_{k,i,\xi}=0\quad{\rm on}\ \{|\xi|=\xi_0\}.$$
				Hence,
				$$(C_{k,i}\psi_k\pm F_{k,i})_{\xi}=0\quad{\rm on}\ \{|\xi|=\xi_0\}.$$
				By the maximum principle and Hopf's lemma, we have
				\begin{equation}\label{Fki maximum principle}
					C_{k,i}\psi_k\pm F_{k,i}\ge\min\left\{0,\min_{\eta=\eta_0\ {\rm or}\ \pi-\delta}C_{k,i}\psi_k\pm F_{k,i}\right\}.
				\end{equation}
				When $\eta=\eta_0$, since $\psi_k(\eta_0)=1$, by the definition of $a_{k,i}$, we have
				\begin{equation}\label{Fki eta0}
					C_{k,i}\psi_k(\eta_0)\pm F_{k,i}(\xi,\eta_0)\ge0.
				\end{equation}
			
				We claim that
				\begin{equation}
					\lim_{\eta\to\pi}C_{k,i}\psi_k\pm F_{k,i}=0,\nonumber
				\end{equation}
				and this is uniform with respect to $\xi$. In fact, since
				$$\int_{\mathbb{S}^{n-2}}u(\xi,\pi,\omega)Y_{k,i}(\omega)\ {\rm d}\omega=0,$$
				we have
				\begin{equation}\label{uniform estimate for Fki}
					\begin{aligned}
						|F_{k,i}(\xi,\eta)|&=\left|\int_{\mathbb{S}^{n-2}}u(\xi,\eta,\omega)Y_{k,i}(\omega)\ {\rm d}\omega\right|\\
						&=\left|\int_{\mathbb{S}^{n-2}}(u(\xi,\eta,\omega)-u(\xi,\pi,\omega))Y_{k,i}(\omega)\ {\rm d}\omega\right|\\
						&\lesssim\Vert Y_{k,i}\Vert_{L^{\infty}}\Vert u_{\eta}\Vert_{L^{\infty}}|\pi-\eta|.
					\end{aligned}
				\end{equation}
				By (\ref{xi eta}), (\ref{bigradient}), Lemma 1.1 and Remark 1.1, we have
				$$|u_{\eta}|^2\lesssim\frac{\rho^2}{\sin^2\eta}|\nabla u|^2=\frac{a^2}{(\cosh\xi-\cos\eta)^2}|\nabla u|^2\lesssim\frac{a^2}{(1-\cos\eta_0)^2}|\nabla u|^2\lesssim\varepsilon^{-2}.$$
				Combining this with (\ref{uniform estimate for Fki}) and (\ref{psik near pi}) proves the claim. By (\ref{Fki maximum principle}), (\ref{Fki eta0}) and the claim, and letting $\delta\to0$, we complete the proof. 
			\end{proof}
		
		\subsection{Estimates of $E_{k,i}$}
			Recall
			$$E_{k,i}\coloneqq\sin^2\eta(F_{k,i,\xi}^2+F_{k,i,\eta}^2)+A_{k,i}F_{k,i}^2.$$
			We have
			\begin{lem}
				Let $A_{k,i}=2\lambda_k+3n$. Then
				$$E_{k,i}\le\left[\frac{100}{\Vert Y_{k,i}\Vert_{L^{\infty}}^2}+\frac{100n\lambda_k}{\Vert\nabla_{\mathbb{S}^{n-2}}Y_{k,i}\Vert_{L^{\infty}}^2}\right]a_{k,i}^2\psi_k^2\quad{\rm in}\ R.$$
			\end{lem}
			\begin{proof}
				It can be seen from a direct calculation that
				\begin{equation}
					\begin{aligned}
						L_k(\sin\eta F_{k,i,\xi})=&2\sin^2\eta\cos\eta F_{k,i,\xi\eta}+(n-3)\sin^3\eta F_{k,i,\xi}\\
						&+(n-2)\sin\eta\left[\left(\cos\eta-\frac{\sin^2\eta}{\cosh\xi-\cos\eta}\right)^2-\sin^2\eta\right]F_{k,i,\xi}\\
						&-(n-2)\frac{\sin^4\eta\sinh\xi}{(\cosh\xi-\cos\eta)^2}F_{k,i,\eta},\nonumber
					\end{aligned}
				\end{equation}
				\begin{equation}
					\begin{aligned}
						L_k(\sin\eta F_{k,i,\eta})=&-2\sin^2\eta\cos\eta F_{k,i,\xi\xi}+(n-3)\sin^3\eta F_{k,i,\eta}\\
						&+(n-2)\frac{\sin^2\eta\sinh\xi}{\cosh\xi-\cos\eta}\left(2\cos\eta-\frac{\sin^2\eta}{\cosh\xi-\cos\eta}\right)F_{k,i,\xi}\\
						&+(n-2)\frac{\sin^3\eta}{\cosh\xi-\cos\eta}\left(2\cos\eta-\frac{\sin^2\eta}{\cosh\xi-\cos\eta}\right)F_{k,i,\eta}.\nonumber
					\end{aligned}
				\end{equation}
				Hence, we have
				\begin{equation}
					\begin{aligned}
						&L_kE_{k,i}\\
						={}&2\sin^2\eta(\sin\eta F_{k,i,\eta\eta}+\cos\eta F_{k,i,\eta})^2+2\sin^2\eta(\sin\eta F_{k,i,\xi\xi}-\cos\eta F_{k,i,\eta})^2\\
						+{}&4\sin^2\eta(\sin\eta F_{k,i,\xi\eta}+\cos\eta F_{k,i,\xi})^2\\
						+{}&\left[\lambda_k+2A_{k,i}-2+2(n-2)\left(\cos\eta-\frac{\sin^2\eta}{\cosh\xi-\cos\eta}\right)^2\right]\sin^2\eta F_{k,i,\xi}^2\\
						+{}&\left[\lambda_k+2A_{k,i}-2+2(n-2)\frac{\sin^2\eta\sinh^2\xi}{(\cosh\xi-\cos\eta)^2}\right]\sin^2\eta F_{k,i,\eta}^2\\
						+{}&4(n-2)\frac{\sin\eta\sinh\xi}{\cosh\xi-\cos\eta}\left(\cos\eta-\frac{\sin^2\eta}{\cosh\xi-\cos\eta}\right)\sin^2\eta F_{k,i,\xi}F_{k,i,\eta}\\
						+{}&A_{k,i}\lambda_kF_{k,i}^2.\nonumber
					\end{aligned}
				\end{equation}
				The first three terms of the right hand side are nonnegative. We now use the fourth and fifth terms to control the sixth term.
				
				We have
				$$1+\cos\eta-\frac{\sin^2\eta}{\cosh\xi-\cos\eta}=\frac{(1+\cos\eta)(\cosh\xi-1)}{\cosh\xi-\cos\eta}\ge0$$
				and
				$$\cos\eta-\frac{\sin^2\eta}{\cosh\xi-\cos\eta}\le\cos\eta\le1.$$
				Hence,
				$$\left|\cos\eta-\frac{\sin^2\eta}{\cosh\xi-\cos\eta}\right|\le1.$$ 
				Moreover, since $\sinh\xi_0=a\sim\varepsilon^{1/2}$ is sufficiently small, for $|\xi|<\xi_0$ we have
				\begin{equation}
					\begin{aligned}
						\frac{\sin^2\eta\sinh^2\xi}{(\cosh\xi-\cos\eta)^2}&=\frac{(1-\cos\eta)(1+\cos\eta)(\cosh\xi-1)(\cosh\xi+1)}{(\cosh\xi-\cos\eta)^2}\\
						&\le\frac{(1+\cos\eta)(\cosh\xi+1)}{4}\le\frac{3}{2}.\nonumber
					\end{aligned}
				\end{equation}
				Then by Young's inequality,
				\begin{equation}
					\begin{aligned}
						&\left|4(n-2)\frac{\sin\eta\sinh\xi}{\cosh\xi-\cos\eta}\left(\cos\eta-\frac{\sin^2\eta}{\cosh\xi-\cos\eta}\right)\sin^2\eta F_{k,i,\xi}F_{k,i,\eta}\right|\\
						\le{}&2(n-2)\sin^2\eta F_{k,i,\xi}^2+3(n-2)\sin^2\eta F_{k,i,\eta}^2.
					\end{aligned}\nonumber
				\end{equation}
				Therefore, we have
				\begin{equation}\label{Lk Eki}
					L_kE_{k,i}\ge\lambda_kE_{k,i}+(2A_{k,i}-3n+4)\sin^2\eta(F_{k,i,\xi}^2+F_{k,i,\eta}^2).
				\end{equation}
				Moreover, by (\ref{equation of psi}) and (\ref{estimate of phik}) we have
				\begin{equation}\label{Lk psik^2}
					\begin{aligned}
						L_k(\psi_k^2)&=\left[\lambda_k+2\phi_k^2-2(n-2)\phi_k\left(1+\cos\eta-\frac{\sin^2\eta}{\cosh\xi-\cos\eta}\right)\right]\psi_k^2\\
						&\le(\lambda_k+2k^2)\psi_k^2.
					\end{aligned}
				\end{equation}
			
				Consider the quantity $E_{k,i}-B_{k,i}\psi_k^2$, where $B_{k,i}$ is a constant to be determined. By (\ref{Lk Eki}), (\ref{Lk psik^2}) and Lemma 4.1 we have
				\begin{equation}
					\begin{aligned}
						&L_k(E_{k,i}-B_{k,i}\psi_k^2)\\
						\ge{}&\lambda_kE_{k,i}+(2A_{k,i}-3n+4)\sin^2\eta(F_{k,i,\xi}^2+F_{k,i,\eta}^2)-B_{k,i}\lambda_k\psi_k^2-2k^2B_{k,i}\psi_k^2\\
						={}&(2A_{k,i}-3n+4+\lambda_k)(E_{k,i}-B_{k,i}\psi_k^2)+(2A_{k,i}-3n+4-2k^2)B_{k,i}\psi_k^2\\
						&-A_{k,i}(2A_{k,i}-3n+4)F_{k,i}^2\\
						\ge{}&(2A_{k,i}-3n+4+\lambda_k)(E_{k,i}-B_{k,i}\psi_k^2)\\
						&+\left[B_{k,i}(2A_{k,i}-3n+4-2k^2)-A_{k,i}(2A_{k,i}-3n+4)\frac{a_{k,i}^2}{\Vert\nabla_{\mathbb{S}^{n-2}}Y_{k,i}\Vert_{L^{\infty}}^2}\right]\psi_k^2.
					\end{aligned}\nonumber
				\end{equation}
				Recall $A_{k,i}=2\lambda_k+3n$ and set
				$$B_{k,i}=\left[\frac{100}{\Vert Y_{k,i}\Vert_{L^{\infty}}^2}+\frac{100n\lambda_k}{\Vert\nabla_{\mathbb{S}^{n-2}}Y_{k,i}\Vert_{L^{\infty}}^2}\right]a_{k,i}^2.$$
				We have
				$$L_k(E_{k,i}-B_{k,i}\psi_k^2)\ge(5\lambda_k+3n+4)(E_{k,i}-B_{k,i}\psi_k^2)\quad{\rm in}\ R_{\delta},$$
				where
				$$R_{\delta}=\{(\xi,\eta):|\xi|<\xi_0,\ \eta_0<\eta<\pi-\delta\}.$$
				By (\ref{Fki}), we have
				$$(E_{k,i}-B_{k,i}\psi_k^2)_{\xi}=0\quad{\rm on}\ \{|\xi|=\xi_0\}.$$
				Hence, the maximum principle and Hopf's lemma imply
				\begin{equation}\label{E-psi maximum principle}
					E_{k,i}-B_{k,i}\psi_k^2\le\max\left\{0,\max_{\eta=\eta_0\ {\rm or}\ \pi-\delta}E_{k,i}-B_{k,i}\psi_k^2\right\}.
				\end{equation}
				When $\eta=\eta_0$, since $\psi_k(\eta_0)=1$, by the definition of $B_{k,i}$, we have
				\begin{equation}\label{Eki eta0}
					E_{k,i}-B_{k,i}\psi_k^2\le0.
				\end{equation}
			
				We claim that
				\begin{equation}
					\lim_{\eta\to\pi}E_{k,i}-B_{k,i}\psi_k^2=0,\nonumber
				\end{equation}
				and this is uniform with respect to $\xi$. In fact, as in the proof of Lemma 4.1, we have
				$$|u_{\xi}|\lesssim\varepsilon^{-1}.$$
				Therefore,
				$$|F_{k,i,\xi}|=\left|\int_{\mathbb{S}^{n-2}}u_{\xi}Y_{k,i}(\omega)\ {\rm d}\omega\right|\lesssim\varepsilon^{-1}\Vert Y_{k,i}\Vert_{L^{\infty}},$$
				and $|F_{k,i,\eta}|$ has the same upper bound. Together with (\ref{psik near pi}) and (\ref{uniform estimate for Fki}), we prove the claim. By (\ref{E-psi maximum principle}), (\ref{Eki eta0}) and the claim, and letting $\delta\to0$, we complete the proof.
			\end{proof}
		
			\subsection{Estimates of $a_{k,i}$}
			As we can see in Lemmas 4.1 and 4.2, we also need to estimate $a_{k,i}$. Recall
			\begin{equation}
				\begin{aligned}
					a_{k,i}\coloneqq\sup_{|\xi|\le\xi_0,\eta=\eta_0}\bigg[\sin\eta_0\Vert Y_{k,i}\Vert_{L^{\infty}}(|F_{k,i,\xi}|+|F_{k,i,\eta}|)+\Vert \nabla_{\mathbb{S}^{n-2}}Y_{k,i}\Vert_{L^{\infty}}|F_{k,i}|\bigg].
				\end{aligned}\nonumber
			\end{equation}
			We have
			\begin{lem}
				$$a_{k,i}\lesssim\left(\Vert Y_{k,i}\Vert_{L^{\infty}}+\Vert \nabla_{\mathbb{S}^{n-2}}Y_{k,i}\Vert_{L^{\infty}}\right)\frac{\Vert\nabla_{\mathbb{S}^{n-2}}Y_{k,i}\Vert_{L^{\infty}}}{\lambda_k^{10n+1}}.$$
			\end{lem}
			\begin{proof}
				Integrating by parts, we have
				\begin{equation}
					\begin{aligned}
						|F_{k,i}(\xi,\eta_0)|&=\left|\int_{\mathbb{S}^{n-2}}u(\xi,\eta_0,\omega)Y_{k,i}\ {\rm d}\omega\right|\\
						&=\frac{1}{\lambda_k}\left|\int_{\mathbb{S}^{n-2}}\Delta_{\mathbb{S}^{n-2}}u(\xi,\eta_0,\omega)Y_{k,i}\ {\rm d}\omega\right|\\
						&=\frac{1}{\lambda_k^{10n+1}}\left|\int_{\mathbb{S}^{n-2}}\nabla_{\mathbb{S}^{n-2}}(\Delta_{\mathbb{S}^{n-2}}^{10n}u(\xi,\eta_0,\omega))\cdot\nabla_{\mathbb{S}^{n-2}}Y_{k,i}\ {\rm d}\omega\right|\\
						&\lesssim\frac{\Vert\nabla_{\mathbb{S}^{n-2}}Y_{k,i}\Vert_{L^{\infty}}}{\lambda_k^{10n+1}}\sup_{\omega\in\mathbb{S}^{n-2}}\left|\nabla_{\mathbb{S}^{n-2}}^{20n+1}u(\xi,\eta_0,\omega)\right|.\nonumber
					\end{aligned}
				\end{equation}
				By the chain rule, for any positive integer $t$, we can check
				$$|\nabla_{\mathbb{S}^{n-2}}^t u|\lesssim\sum_{j=1}^t|\nabla^j u(x)|.$$
				By the argument in Section 2, the domain corresponding to $\{\eta=\eta_0\}$ in the Euclidean coordinates is contained in $\Omega\setminus\overline{\Omega_0\cup B_1\cup B_2}$, so we can use Lemma 1.1 to get
				\begin{equation}
					|F_{k,i}(\xi,\eta_0)|\lesssim\frac{\Vert\nabla_{\mathbb{S}^{n-2}}Y_{k,i}\Vert_{L^{\infty}}}{\lambda_k^{10n+1}}.\nonumber
				\end{equation}
				Similarly, we have
				$$|F_{k,i,\xi}(\xi,\eta_0)|\lesssim\frac{\Vert\nabla_{\mathbb{S}^{n-2}}Y_{k,i}\Vert_{L^{\infty}}}{\lambda_k^{10n+1}}\sup_{\omega\in\mathbb{S}^{n-2}}\left|\nabla_{\mathbb{S}^{n-2}}^{20n+1}u_{\xi}(\xi,\eta_0,\omega)\right|.$$
				By the chain rule,
				$$u_{\xi}=u_{\rho}\rho_{\xi}+u_{x_n}x_{n,\xi}.$$
				It can be seen from (\ref{xi eta}) that
				\begin{equation}\label{estimate of a/D}
					|\rho_{\xi}|^2+|x_{n,\xi}|^2=\frac{a^2}{(\cosh\xi-\cos\eta)^2}\le\frac{a^2}{(1-\cos\eta_0)^2}\lesssim\varepsilon^{-1}.
				\end{equation}
				Using the chain rule once again, for any positive integer $t$, we have
				$$|\nabla_{\mathbb{S}^{n-2}}^t u_{\xi}(\xi,\eta,\omega)|\lesssim\varepsilon^{-1/2}\sum_{j=1}^{t+1}|\nabla^j u(x)|.$$
				Hence,
				$$|F_{k,i,\xi}(\xi,\eta_0)|\lesssim\varepsilon^{-1/2}\frac{\Vert\nabla_{\mathbb{S}^{n-2}}Y_{k,i}\Vert_{L^{\infty}}}{\lambda_k^{10n+1}}.$$
				The estimate of $|F_{k,i,\eta}(\xi,\eta_0)|$ can be proved in the same way. Since
				$$\sin\eta_0=\frac{8ar_0}{a^2+16r_0^2}\sim\varepsilon^{1/2},$$
				the proof is complete.
			\end{proof}
			
			Now we can give an estimate for $|\nabla u|$.
			\begin{prop}
				$$|\nabla u|\lesssim|\nabla u_0|+(\varepsilon+\rho^2)^{\frac{\alpha_1-1}{2}}\quad{\rm in}\ \Omega_0.$$
			\end{prop}
			\begin{proof}
				Recall
				$$A_{k,i}=2\lambda_k+3n,\quad B_{k,i}=\left[\frac{100}{\Vert Y_{k,i}\Vert_{L^{\infty}}^2}+\frac{100n\lambda_k}{\Vert\nabla_{\mathbb{S}^{n-2}}Y_{k,i}\Vert_{L^{\infty}}^2}\right]a_{k,i}^2.$$
				We can see from Lemma 4.3 that
				\begin{equation}
					\begin{aligned}
						B_{k,i}\lesssim\frac{1}{\lambda_k^{20n+2}}\bigg[&\Vert\nabla_{\mathbb{S}^{n-2}}Y_{k,i}\Vert_{L^{\infty}}^2+\lambda_k\Vert Y_{k,i}\Vert_{L^{\infty}}^2\\
						&+\lambda_k\Vert\nabla_{\mathbb{S}^{n-2}}Y_{k,i}\Vert_{L^{\infty}}^2+\frac{\Vert\nabla_{\mathbb{S}^{n-2}}Y_{k,i}\Vert_{L^{\infty}}^4}{\Vert Y_{k,i}\Vert_{L^{\infty}}^2}\bigg].
					\end{aligned}\nonumber
				\end{equation}
				Moreover, we have
				$$\lambda_k\sim k^2,\quad1\lesssim\Vert Y_{k,i}\Vert_{L^{\infty}}\lesssim k^{\frac{n-3}{2}},\quad\Vert\nabla_{\mathbb{S}^{n-2}}Y_{k,i}\Vert_{L^{\infty}}\lesssim k^{\frac{n-1}{2}}.$$
				These estimates are trivial when $n=3$ and follow from \cite{estimate Yki} for $n\ge4$. Therefore,
				\begin{equation}\label{estimates of B_{k,i}}
					B_{k,i}\lesssim k^{-16n}.
				\end{equation}
				Then (\ref{bigradient}), (\ref{estimates of B_{k,i}}) and Lemma 4.2 imply
				\begin{equation}\label{simple estimates of nabla uki}
					\begin{aligned}
						|\nabla u_{k,i}|&\lesssim\frac{1}{\rho}\left(\Vert Y_{k,i}\Vert_{L^{\infty}}+\frac{\Vert \nabla_{\mathbb{S}^{n-2}}Y_{k,i}\Vert_{L^{\infty}}}{A_{k,i}^{1/2}}\right)E_{k,i}^{1/2}\\
						&\lesssim\frac{1}{\rho}\left(\Vert Y_{k,i}\Vert_{L^{\infty}}+\frac{\Vert \nabla_{\mathbb{S}^{n-2}}Y_{k,i}\Vert_{L^{\infty}}}{A_{k,i}^{1/2}}\right)B_{k,i}^{1/2}\psi_k\\
						&\lesssim\rho^{-1}k^{-5n}\psi_k.
					\end{aligned}
				\end{equation}
				So it suffices to estimate $\psi_k$.
				
				We can see from Lemma 3.4 that
				\begin{equation}
					\psi_k\le\Psi_k=\left(\frac{a\cot\frac{\eta}{2}}{4r_0}\right)^k\left(\frac{a^2+16r_0^2}{a^2+a^2\cot^2\frac{\eta}{2}}\right)^{\frac{k-\alpha_k}{2}}.\nonumber
				\end{equation}
				For $k=1$, we can check $\alpha_1\in(0,1)$. By (\ref{a}) and (\ref{rho eta}), we have
				$$a\sim\varepsilon^{1/2},\quad\rho\sim a\cot\frac{\eta}{2}\quad{\rm in}\ R.$$
				Hence,
				$$|\nabla u_{1,i}|\lesssim\rho^{-1}\psi_1\lesssim(\varepsilon+\rho^2)^{\frac{\alpha_1-1}{2}}.$$
				For $k\ge2$, we can check $k-\alpha_k\in(0,1)$. Since $$a\cot\frac{\eta_0}{2}=4r_0,$$
				we have
				\begin{equation}
					\begin{aligned}
						\psi_k&\lesssim\rho\left(\frac{a\cot\frac{\eta}{2}}{4r_0}\right)^{k-1}\left(\frac{a^2+16r_0^2}{a^2+a^2\cot^2\frac{\eta}{2}}\right)^{\frac{k-\alpha_k}{2}}\\
						&=\rho\left(\frac{a\cot\frac{\eta}{2}}{4r_0}\right)^{\alpha_k-1}\left(1+\frac{a^2}{16r_0^2}\right)^{\frac{k-\alpha_k}{2}}\left(\frac{a^2\cot^2\frac{\eta}{2}}{a^2+a^2\cot^2\frac{\eta}{2}}\right)^{\frac{k-\alpha_k}{2}}\\
						&\lesssim\rho\left(1+\frac{a^2}{16r_0^2}\right)^{\frac{1}{2}}\lesssim\rho.
					\end{aligned}\nonumber
				\end{equation}
				Hence,
				$$|\nabla u_{k,i}|\lesssim k^{-5n}.$$
				Since $N(k)\sim k^{n-3}$, by the above argument we have
				\begin{equation}
					\begin{aligned}
						|\nabla u|&\lesssim|\nabla u_0|+\sum_{k=1}^{\infty}\sum_{i=1}^{N(k)}|\nabla u_{k,i}|\\
						&\lesssim|\nabla u_0|+(\varepsilon+\rho^2)^{\frac{\alpha_1-1}{2}}+\sum_{k=2}^{\infty}\frac{N(k)}{k^{5n}}\\
						&\lesssim|\nabla u_0|+(\varepsilon+\rho^2)^{\frac{\alpha_1-1}{2}}.
					\end{aligned}\nonumber
				\end{equation}
			\end{proof}

	\section{Estimates of $|\nabla u_0|$}
		In this section, we return to the cylindrical coordinates and show that $|\nabla u_0|$ is uniformly bounded.
		
		Define
		$$\tilde{B}_i\coloneqq\{(\rho,x_n):\rho^2+\left[x_n+(-1)^i(1+\varepsilon)\right]^2<1\}$$
		and
		$$\Gamma_i=\{(\rho,x_n):0\le\rho\le r_0,\ x_n=(-1)^{i+1}h(\rho)\}.$$
		We point out that in the $\rho Ox_n$ half-plane, the domain $\tilde{B}_i$ represents a semidisk, while $\Gamma_i$ refers to a part of the arc in the narrow region. Recall
		$$h(\rho)\coloneqq1+\varepsilon-\sqrt{1-\rho^2}$$
		and define
		$$\tilde{\Omega}_0\coloneqq\{(\rho,x_n):0\le\rho<r_0,\ |x_n|<h(\rho)\}.$$
		
		Since $u_0=u_0(\rho,x_n)$, we consider the problem in the $\rho Ox_n$ half-plane. In the cylindrical coordinates, we have
		\begin{equation}\label{cylgradient}
			|\nabla u_0|^2=|u_{0,\rho}|^2+|u_{0,x_n}|^2.
		\end{equation}
		Hence, it suffices to estimate $|u_{0,\rho}|$ and $|u_{0,x_n}|$.
		
		By the decomposition (\ref{fourier series}), $u_0$ can be seen as a function of $\rho$ and $x_n$:
		$$u_0(\rho,x_n)=\fint_{\mathbb{S}^{n-2}}u(\rho,\omega,x_n)\ {\rm d}\omega.$$
		Moreover, it is harmonic with respect to $x=(x',x_n)$. Therefore,
		\begin{equation}\label{equation of u0}
			\begin{cases}
				\displaystyle u_{0,\rho\rho}+\frac{n-2}{\rho}u_{0,\rho}+u_{0,x_nx_n}=0\quad{\rm in}\ \tilde{\Omega}_0,\\
				\partial_{\nu}u_0=0\quad{\rm on}\ \Gamma_1\cup\Gamma_2,
			\end{cases}
		\end{equation}
		where $\nu$ denotes the normal vector pointing into $\tilde{B}_i$. It can be seen from (\ref{equation of u0}) that
		$$(\rho^{n-2}u_{0,\rho})_{\rho}+(\rho^{n-2}u_{0,x_n})_{x_n}=0.$$
		Hence, the following integral is independent of the path:
		\begin{equation}\label{U0}
			U_0(\rho,x_n)=\int_{(0,0)}^{(\rho,x_n)}-\tilde{\rho}^{n-2}u_{0,x_n}\ {\rm d}\tilde{\rho}+\tilde{\rho}^{n-2}u_{0,\rho}\ {\rm d}\tilde{x}_n.
		\end{equation}
		A direct calculation shows
		\begin{equation}\label{derivative of U}
			U_{0,\rho}=-\rho^{n-2}u_{0,x_n},\quad U_{0,x_n}=\rho^{n-2}u_{0,\rho},
		\end{equation}
		and
		\begin{equation}\label{equation of U}
			U_{0,\rho\rho}-\frac{n-2}{\rho}U_{0,\rho}+U_{0,x_nx_n}=0.
		\end{equation}
		
		We first use the maximum principle to give an estimate of $U_0$. Then we obtain estimates for $|\nabla U_0|$ and $|\nabla u_0|$.
		\begin{lem}
			Let
			$$V_0(\rho,x_n)\coloneqq\rho^{n-1}[h^2(\rho)-x_n^2].$$
			Then
			$$|U_0|\lesssim V_0\quad{\rm in}\ \tilde{\Omega}_0.$$
		\end{lem}
		\begin{proof}
			Let $\tau$ be the unit tangent vector of $\Gamma_1\cup\Gamma_2$. By (\ref{derivative of U}) and the boundary condition in (\ref{equation of u0}), we have
			$$\partial_{\tau}U_0=0\quad{\rm on}\ \Gamma_1\cup\Gamma_2.$$
			Hence, $U_0$ is constant on $\Gamma_i$. Moreover, we see from (\ref{U0}) that $U_0(0,x_n)\equiv0$. Therefore,
			$$U_0\equiv0\quad{\rm on}\ \Gamma_1\cup\Gamma_2\cup\{(\rho,x_n):\rho=0,\ |x_n|\le\varepsilon\},$$
			and it is easy to check that $V_0$ also satisfies this condition.
			
			Since
			$$h'(\rho)=\frac{\rho}{(1-\rho^2)^{1/2}},\quad h''(\rho)=\frac{1}{(1-\rho^2)^{3/2}},$$
			we have
			\begin{equation}\label{equation of V}
				\begin{aligned}
					&V_{0,\rho\rho}-\frac{n-2}{\rho}V_{0,\rho}+V_{0,x_nx_n}\\
					={}&\rho^{n-1}\left[2n\rho^{-1}hh'+2(h')^2+2hh''-2\right]\\
					\le{}&-\rho^{n-1}\le0.
				\end{aligned}
			\end{equation}
			Therefore,
			$$(C_0V_0\pm U_0)_{\rho\rho}-\frac{n-2}{\rho}(C_0V_0\pm U_0)_{\rho}+(C_0V_0\pm U_0)_{x_nx_n}\le0\quad{\rm in}\ \tilde{\Omega}_{\delta},$$
			where
			$$C_0=\frac{\Vert u_{0,\rho}\Vert_{L^{\infty}(\{\rho=r_0\})}}{r_0h(r_0)}$$
			and
			$$\tilde{\Omega}_{\delta}=\{(\rho,x_n):\delta<\rho<r_0,\ |x_n|<h(\rho)\}.$$
			We point out that
			$$|u_{0,\rho}(r_0,x_n)|\lesssim\left|\fint_{\mathbb{S}^{n-2}}u_{\rho}(r_0,\omega,x_n)\ {\rm d}\omega\right|\lesssim\Vert \nabla u\Vert_{L^{\infty}(\{|x'|=r_0\})}\lesssim1,$$
			which means $|C_0|\lesssim1$.
			
			In $\tilde{\Omega}_{\delta}$, by the maximum principle we have
			\begin{equation}\label{DV-U maximum principle}
				C_0V_0\pm U_0\ge\min_{\partial\tilde{\Omega}_{\delta}}C_0V_0\pm U_0=\min\{0,\min_{\rho=\delta\ {\rm or}\ r_0}C_0V_0\pm U_0\}.
			\end{equation}
			When $\rho=r_0$, if $x_n\ge0$, by (\ref{derivative of U}) we have
			\begin{equation}\label{DV-U near r0}
				\begin{aligned}
					|U_0(r_0,x_n)|&=|U_0(r_0,x_n)-U_0(r_0,h(r_0))|\\
					&\le r_0^{n-2}\Vert u_{0,\rho}\Vert_{L^{\infty}(\{\rho=r_0\})}|h(r_0)-|x_n||\le C_0V_0(r_0,x_n).
				\end{aligned}
			\end{equation}
			If $x_n<0$, by considering $U_0(r_0,x_n)-U_0(r_0,-h(r_0))$ we can get the same estimate.
			 
			We claim that
			$$\lim_{\rho\to0}C_0V_0\pm U_0=0,$$
			and this is uniform with respect to $x_n$. In fact, we have $$|V_0|\lesssim\rho^{n-1}$$
			and
			\begin{equation}
				\begin{aligned}
					|U_0(\rho,x_n)|&=|U_0(\rho,x_n)-U_0(\rho,h(\rho))|\\
					&=\left|\int_{x_n}^{h(\rho)}U_{0,x_n}(\rho,t)\ {\rm d}t\right|\le2h(\rho)\rho^{n-2}\Vert u_{0,\rho}\Vert_{L^{\infty}(\Omega_0)}.
				\end{aligned}\nonumber
			\end{equation}
			Since
			$$|u_{0,\rho}|=\left|\fint_{\mathbb{S}^{n-2}}u_{\rho}(\rho,\omega,x_n)\ {\rm d}\omega\right|\lesssim\Vert \nabla u\Vert_{L^{\infty}(\Omega_0)}\lesssim\varepsilon^{-1/2},$$
			the claim holds. By (\ref{DV-U maximum principle}), (\ref{DV-U near r0}) and the claim, and letting $\delta\to0$, we complete the proof.
		\end{proof}
		Now we can give an estimate of $|\nabla u_0|$.
		\begin{prop}
			$$|\nabla u_0|\lesssim 1\quad{\rm in}\ \Omega_0.$$
		\end{prop}
		\begin{proof}
			Since
			$$U_0=V_0=0\quad{\rm on}\ \Gamma_1\cup\Gamma_2,$$
			Lemma 5.1 implies that on $\ \Gamma_1\cup\Gamma_2$, we have
			$$|U_{0,\rho}|^2+|U_{0,x_n}|^2\lesssim|\partial_{\nu}U_0|^2\lesssim|\partial_{\nu}V_0|^2\lesssim|V_{0,\rho}|^2+|V_{0,x_n}|^2\lesssim\rho^{2n-4}.$$
			Now we can see from (\ref{cylgradient}) and (\ref{derivative of U}) that
			$$|\nabla u_0|^2\lesssim\rho^{4-2n}(|U_{0,\rho}|^2+|U_{0,x_n}|^2)\lesssim1\quad{\rm on}\ \partial\Omega_0\cap(\partial B_1\cup\partial B_2).$$
			Moreover, by Lemma 1.1 we have
			$$|\nabla u_0|\lesssim1\quad{\rm on}\ \partial\Omega_0\setminus(\partial B_1\cup\partial B_2).$$
			Since $|\nabla u_0|^2$ is subharmonic, by the maximum principle, $|\nabla u_0|$ is uniformly bounded in $\Omega_0$.
		\end{proof}	
		Finally, note that the constant $\alpha_1$ is the same as $\alpha$ in Theorem 1.1. By Propositions 4.1 and 5.1, we prove Theorem 1.1.

		\section*{Acknowledgements}
		The author acknowledges the use of OpenAI tools in the preparation of this manuscript. The bispherical coordinates were initially proposed with the aid of ChatGPT. Then the author prompted AI to consider the eigenvalue problem on $\mathbb{S}^{n-2}$. AI gave a proof for $n=3$ but ran into difficulties in higher dimensions. The author learnt its argument and found that the construction of $\psi_k$ and the estimate of $|\nabla u_0|$ can be extended to higher dimensions easily. The author's main work is in Section 4 and the latter part of Section 3. The author refined the method given by AI and completed the proof in all dimensions. The author assumes full responsibility for the mathematical validity, accuracy, and integrity of the proofs presented in the manuscript.


\begin{thebibliography}{30}
		\bibitem{Ammariqqb}
		H.~Ammari, G.~Ciraolo, H.~Kang, H.~Lee, and K.~Yun,
		Spectral analysis of the Neumann--Poincar\'e operator and characterization
		of the stress concentration in anti-plane elasticity,
		\emph{Arch. Ration. Mech. Anal.} \textbf{208} (2013), 275--304.
		
		\bibitem{Ammari2}
		H.~Ammari, H.~Kang, H.~Lee, J.~Lee, and M.~Lim,
		Optimal estimates for the electric field in two dimensions,
		\emph{J. Math. Pures Appl. (9)} \textbf{88} (2007), 307--324.
		
		\bibitem{Ammari1}
		H.~Ammari, H.~Kang, and M.~Lim,
		Gradient estimates for solutions to the conductivity problem,
		\emph{Math. Ann.} \textbf{332} (2005), 277--286.
		
		\bibitem{perfect-alldim}
		E.~S.~Bao, Y.~Y.~Li, and B.~Yin,
		Gradient estimates for the perfect conductivity problem,
		\emph{Arch. Ration. Mech. Anal.} \textbf{193} (2009), 195--226.
		
		\bibitem{insulated-dim2}
		E.~S.~Bao, Y.~Y.~Li, and B.~Yin,
		Gradient estimates for the perfect and insulated conductivity problems
		with multiple inclusions,
		\emph{Comm. Partial Differential Equations} \textbf{35} (2010), 1982--2006.
		
		\bibitem{green function}
		H.~Dong and H.~Li,
		Optimal estimates for the conductivity problem by Green's function method,
		\emph{Arch. Ration. Mech. Anal.} \textbf{231} (2019), 1427--1453.
		
		\bibitem{insulated-alldim general}
		H.~Dong, Y.~Y.~Li, and Z.~Yang,
		Gradient estimates for the insulated conductivity problem:
		the non-umbilical case,
		\emph{J. Math. Pures Appl. (9)} \textbf{189} (2024),
		Paper No.~103587, 37~pp.
		
		\bibitem{insulated-alldim}
		H.~Dong, Y.~Y.~Li, and Z.~Yang,
		Optimal gradient estimates of solutions to the insulated conductivity
		problem in dimension greater than two,
		\emph{J. Eur. Math. Soc. (JEMS)} \textbf{27} (2025), 3275--3296.
		
		\bibitem{KLYP}
		H.~Kang, M.~Lim, and K.~Yun,
		Asymptotics and computation of the solution to the conductivity equation
		in the presence of adjacent inclusions with extreme conductivities,
		\emph{J. Math. Pures Appl. (9)} \textbf{99} (2013), 234--249.
		
		\bibitem{generalsys}
		Y.~Y.~Li and L.~Nirenberg,
		Estimates for elliptic systems from composite material,
		\emph{Comm. Pure Appl. Math.} \textbf{56} (2003), 892--925.
		
		\bibitem{generalk} 
		Y.~Y.~Li and M.~Vogelius,
		Gradient estimates for solutions to divergence form elliptic equations
		with discontinuous coefficients,
		\emph{Arch. Ration. Mech. Anal.} \textbf{153} (2000), 91--151.
		
		\bibitem{maximum principle with general k}
		Y.~Y.~Li and B.~Weinkove,
		Composite media, almost touching disks and the maximum principle,
		\emph{Pure Appl. Math. Q.} \textbf{21} (2025), 2287--2301.
		
		\bibitem{insulated-Harnack}
		Y.~Y.~Li and Z.~Yang,
		Gradient estimates of solutions to the insulated conductivity problem
		in dimension greater than two,
		\emph{Math. Ann.} \textbf{385} (2023), 1775--1796.
		
		\bibitem{ode lemma}
		G.~Teschl,
		\emph{Ordinary Differential Equations and Dynamical Systems},
		Graduate Studies in Mathematics, vol.~140,
		American Mathematical Society, Providence, RI, 2012.
		
		\bibitem{maxprin}
		B.~Weinkove,
		The insulated conductivity problem, effective gradient estimates
		and the maximum principle,
		\emph{Math. Ann.} \textbf{385} (2023), 1--16.
		
		\bibitem{estimate Yki}
		B.~Xu,
		Derivatives of the spectral function and Sobolev norms of
		eigenfunctions on a closed Riemannian manifold,
		\emph{Ann. Global Anal. Geom.} \textbf{26} (2004), 231--252.
		
		\bibitem{perfect-dim2-general}
		K.~Yun,
		Estimates for electric fields blown up between closely adjacent
		conductors with arbitrary shape,
		\emph{SIAM J. Appl. Math.} \textbf{67} (2007), 714--730.
		
		\bibitem{Yun another}
		K.~Yun,
		Optimal bound on high stresses occurring between stiff fibers
		with arbitrary shaped cross-sections,
		\emph{J. Math. Anal. Appl.} \textbf{350} (2009), 306--312.
		
		\bibitem{Yun Ball}
		K.~Yun,
		An optimal estimate for electric fields on the shortest line segment
		between two spherical insulators in three dimensions,
		\emph{J. Differential Equations} \textbf{261} (2016), 148--188.
		
	\end{thebibliography}
\end{document}